\documentclass[12pt]{amsart}
\usepackage{array,amsmath, enumerate, color, url}
\usepackage{amssymb}
\usepackage{graphicx,subfigure}

\makeatletter
 \usepackage{fullpage} % reduces the margins to 1 inch
 \usepackage{amsthm}
\usepackage{fullpage}

\makeatother

\begin{document}

\newtheorem{thm}{Theorem}[section]
\newtheorem{lem}[thm]{Lemma}
\newtheorem{prop}[thm]{Proposition}
\newtheorem{cor}[thm]{Corollary}
\newtheorem{con}[thm]{Conjecture}
\newtheorem{claim}[thm]{Claim}
\newtheorem{obs}[thm]{Observation}
\newtheorem{definition}[thm]{Definition}
\newtheorem{example}[thm]{Example}
\newtheorem{rmk}[thm]{Remark}
\newcommand{\di}{\displaystyle}
\def\dfc{\mathrm{def}}
\def\cF{{\cal F}}
\def\cH{{\cal H}}
\def\cT{{\cal T}}
\def\C{{\mathcal C}}
\def\cA{{\cal A}}
\def\cB{{\mathcal B}}
\def\P{{\mathcal P}}
\def\Q{{\mathcal Q}}
\def\cP{{\mathcal P}}
\def\cp{\alpha'}
\def\Frk{F_k^{2r+1}}

\title{Optimal open-locating-dominating sets in infinite triangular grids}
\author{Rex Kincaid \and Allison Oldham \and Gexin Yu}
\thanks{The third author's research was supported in part by NSA grant H98230-12-1-0226. }
\address{Department of Mathematics, The College of William and Mary,
Williamsburg, VA, 23185.}
\email{gyu@wm.edu}

\date{\today}

\maketitle

\begin{abstract}
An open-locating-dominating set (OLD-set) is a subset of vertices of a graph such that every vertex in the 
graph has at least one neighbor in the set and no two vertices in the graph have the same set of 
neighbors in the set.  This is an analogue to the well-studied identifying code in the literature.  
In this paper, we prove that the optimal density of the OLD-set for the infinite triangular grid is $4/13$.
\end{abstract}

\section{Introduction}

Various types of protection sets are motivated by the desire to detect failures efficiently in settings 
where the underlying structure can be modeled as a network or graph.   One such example is location 
detection in a wireless sensor network where one tries to find a distinguishing set. A collection of 
subsets $\{S_1, S_2, \ldots, S_k\}$ of $V(G)$ is a {\em distinguishing set} if $\cup_{i=1}^k S_i=V(G)$
and for every pair of distinct vertices $u, v\in V(G)$, some $S_i$ contains exactly one of them. 
In \cite{RSTU04}, {\em identifying codes} are shown to provide solutions for optimal sensor placement.  
{\em Identifying codes}, or {\em ID-codes},  introduced in 1998 by Karpovsky, Chakrabarty and Levitin 
\cite{KCL98}, is a subset of vertices $S=\{v_1, v_2, \ldots, v_k\}$ such that 
$\{N[v_1], N[v_2], \ldots, N[v_k]\}$ is a distinguishing set, where $N[v]$ is closed neighborhood of 
$v$  (the union of $v$ and its neighborhood $N(v)$).   Seo and Slater \cite{SS10, SS11} introduced 
and studied {\em Open Locating Dominating sets} or {\it OLD-sets}, a subset of vertices 
$S=\{v_1, v_2, \ldots, v_k\}$ such that $\{N(v_1), N(v_2), \ldots, N(v_k)\}$ is a distinguishing set.   
OLD-sets can be used to identify 
malfunctioning elements in a network \cite{Slater13}.    Honkala, Laihonen and Ranto~\cite{HLR02} 
introduced and studied {\em strong identifying codes}, a subset of vertices that is both an 
identifying code and an OLD-set.  

We study the OLD-sets here. As mentioned above, a set $S\subseteq V(G)$ is an OLD-set if

\begin{itemize}
\item For every vertex $v\in V(G)$, $S \cap N(v) \neq \emptyset$.
\item For any vertices $v,w\in V(G)$ such that $v \neq w$, $(N(v) \cap S) \neq (N(w) \cap S)$.
\end{itemize}

Not every graph has an OLD-set or an ID-code. It was shown in \cite{SS10} that an OLD-set exists 
if and only if the graph has no isolated vertices and no two vertices have the same neighborhood.  

For finite graphs, the {\em OLD-number} is defined to be the minimum cardinality of an OLD-set. 
It is shown in \cite{SS10} that it is NP-complete to determine if the OLD-number is at most $k$ 
for a given positive integer $k$.   For infinite graphs, a natural extension of the OLD-number is
the {\em OLD-density} or simply the {\em density} when the context is clear.

\vskip .1 true in
\begin{definition}
The {\em OLD-density} of an OLD-set is defined to be
$$\min_{v \in V(G)} \{ \lim\sup \frac{|N_k[v] \cap S|}{|N_k[v]|}:  \text{ $S$ is an OLD-set}\},$$
where $N_k[v]$ is the set of all vertices in $G$ that are within $k$ distance of $v$.
\end{definition}
\vskip .1 true in

%we define the {\em density} of OLD-sets to be
%$$\min_{v \in V(G)} \{ \lim_{k\to\infty} \frac{|N_k[v] \cap S|}{|N_k\{v\}|}:  \text{ $S$ is an OLD-set}\},$$
%where $N_k[v]$ is the set of all vertices in $G$ that are within $k$ distance of $v$.  

We can think of the {\em OLD-density} as the percentage of vertices of a graph which are included in an OLD-set. 

\begin{figure}[h]
\begin{center}
\includegraphics[scale=0.55]{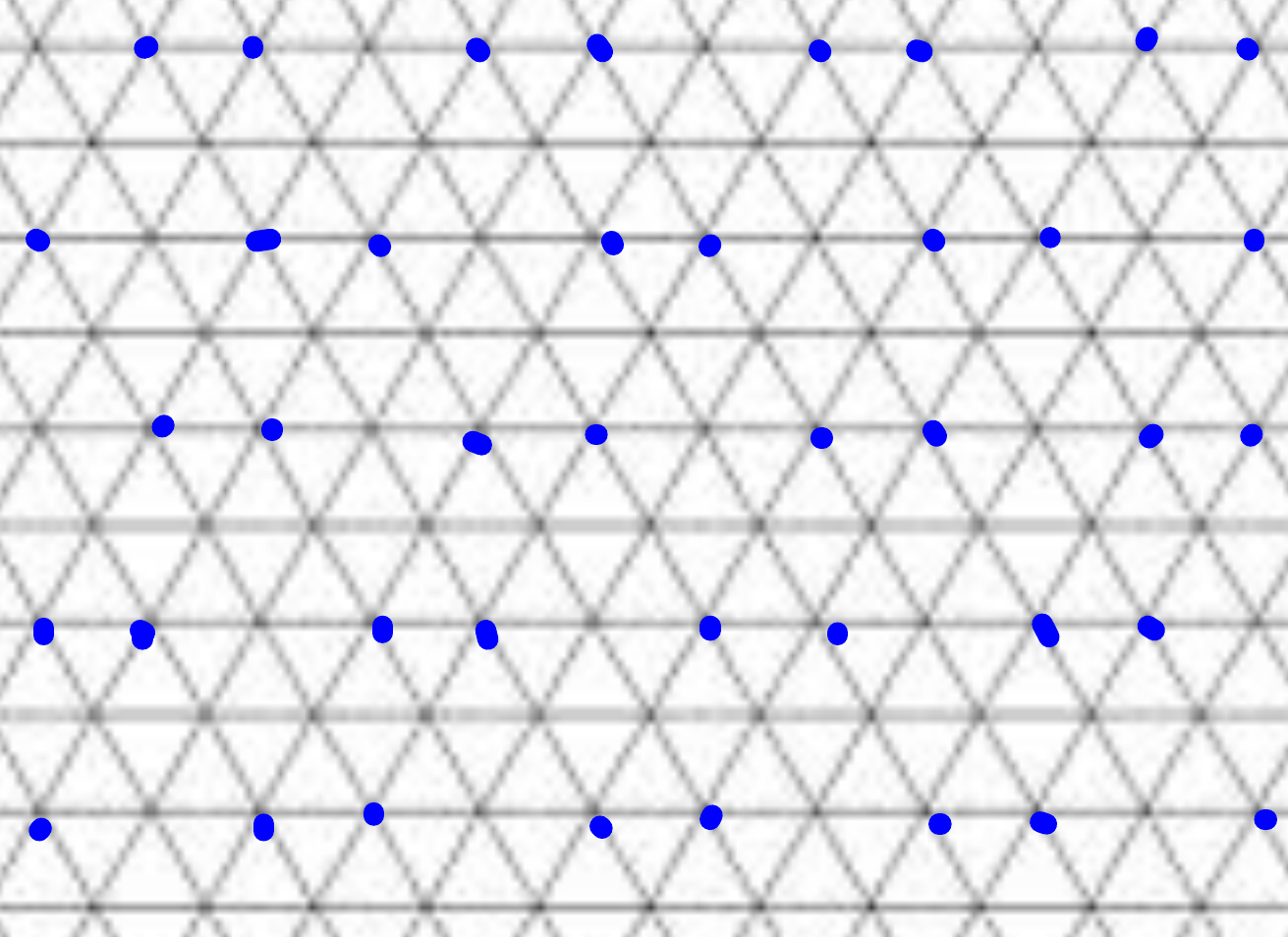}\hskip 0.7in
\includegraphics[scale=0.55]{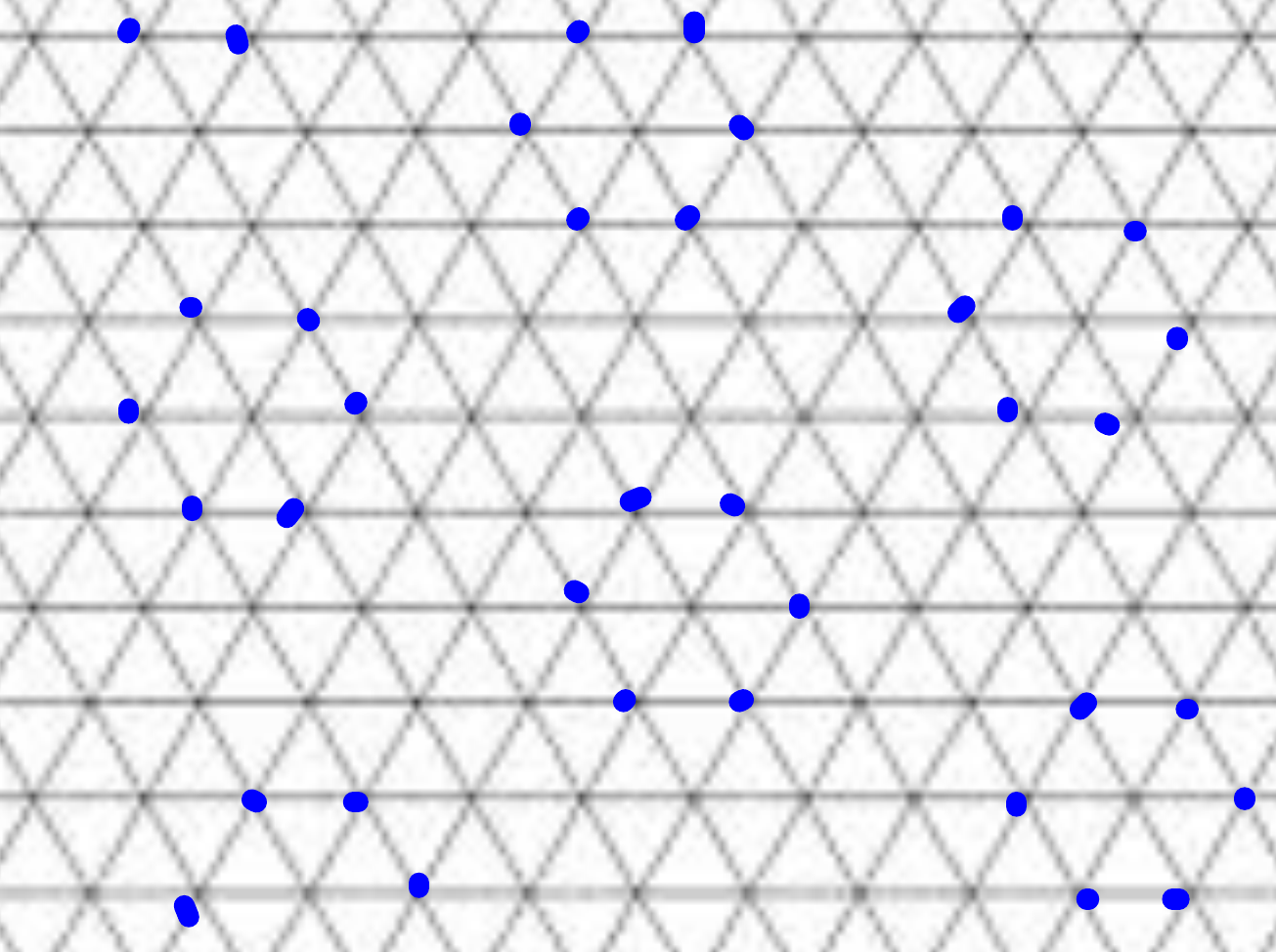}
\caption{OLD-set with density $1/3$ (left) and $6/19$ (right)}
\label{3.1}
\end{center}
\end{figure}

It has been a popular topic to study the density of identifying codes in infinite grids, see 
\cite{AY13, BL05, CGHLM99, CHLZ00, CY09, H09, MS10}.  
It may be difficult to determine the exact density of identifying codes for many infinite grids, for example, 
the infinite hexagonal grid.   
% ## I am not sure what you are trying to say with the above sentence. Do you mean....
% ## "It may be hard to determine the exact density of identifying codes for many infinite grids, for
% ## example, the infinite hexagonal grid."
% ## By hard, do you mean NP-Hard, or that the exact density values have thus far eluded researchers
% ## efforts to find and prove them.
We are also interested in the OLD-density of infinite grids.  
It is not hard to show (see \cite{SS10}) that every $r$-regular graph has OLD-density at least $2/(r+1)$.  
Seo and Slater \cite{SS10} showed that the OLD-density of an infinite square grid is $2/5$ and an infinite hexagonal grid is $1/2$.    
For an infinite triangular grid, they gave a construction with density $1/3$.   Honkala\cite{H09} gave a better construction with a density of $6/19$, 
see Figure~\ref{3.1}.  For more results on identifying codes of OLD-sets, see the dynamic survey by Lobstein~\cite{L}, which lists more than 270 articles on those topics.

\begin{figure}[h]
\begin{center}
\includegraphics[scale=0.8]{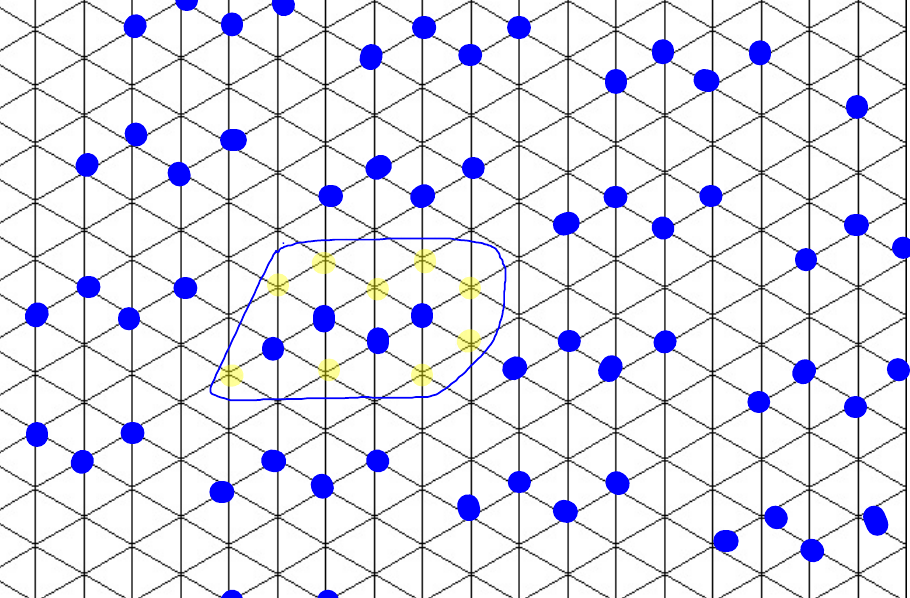}
\caption{A construction with Density $4/13$}
\label{4-13}
\end{center}
\end{figure}

As the infinite triangular grid is $6$-regular, an obvious lower bound for the density of an OLD-set is $\frac{2}{7}$.  
Here we give a construction with density $4/13$, see Figure~\ref{4-13}. 

We also show that the density is at least $4/13$.  Thus we obtain the optimal density of the infinite triangular grid. 

\begin{thm}
The optimal OLD-density of the infinite triangular grid is $4/13$.  
\end{thm}

The proof of the lower bound uses a discharging argument, a popular method in graph coloring problems and recently adapted 
to the study of identifying codes, see \cite{CY09, AY13}.   The idea is to assign weight $1$ to vertices in the OLD-set and 
$0$ to other vertices, and then redistribute the weights among vertices so that each vertex ends up with a weight at least $4/13$.  
This proof is given in the next seOLD-setction.

\section{Proof of the lower bound}

In this section, we show that each OLD-set in an infinite triangular grid has density at least $4/13$.   

Let $S$ be an OLD-set in the infinite triangular grid.  As mentioned above, we assign weight $1$ to each vertex in $S$ and $0$ 
to vertices not in $S$.  We redistribute the weights so that each vertex on the grid has weight at least $4/13$.  Before proceeding
we provide a collection of definitions needed for the proof.

\begin{figure}[h]
\begin{center}
\includegraphics[scale=0.8]{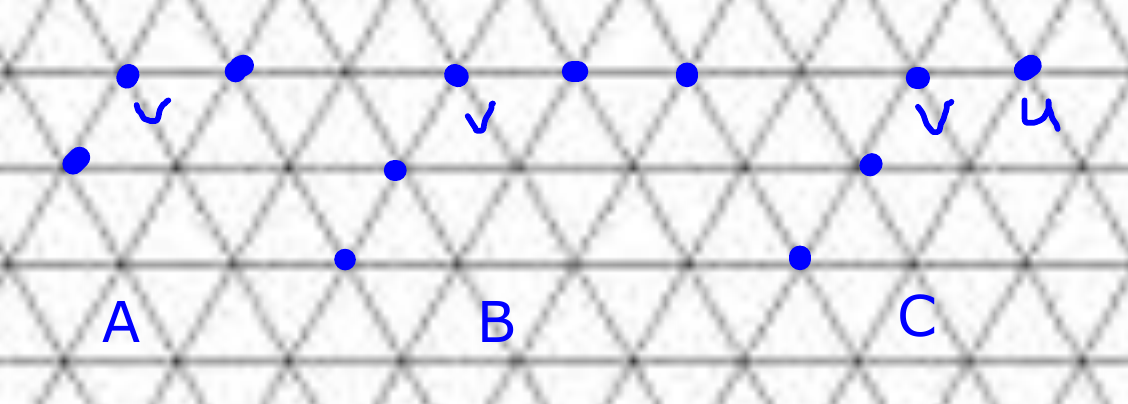}
\caption{Corner $2$-vertex (A), poor $2$-vertex (B), and poor couple (C)}
\label{poor-vertices}
\end{center}
\end{figure}

We define a {\em $t$-cluster} in $S$ to be a connected component in $S$ with $t$ vertices.  Then $t\ge 2$, as $1$-cluster is not  
allowed in an OLD-set.  If $v$ is in cluster $C$, then we denote $d_C(v)$ the number of neighbors $v$ has in $C$ and call $v$ a  
$k$-vertex  (or $k^+$-vertex) if $d_C(v)=k$ (or $d_C(v)\ge k$).  A {\em corner $2$-vertex} in a cluster $C$ is a $2$-vertex whose 
two neighbors in $C$ share two neighbors  in $G$,  and a {\em poor $2$-vertex} is a corner $2$-vertex that has no $3^+$-vertices 
or corner $2$-vertices as neighbors;   a {\em $1,2$-couple} is a pair of adjacent vertices $u,v\in C$ such that $d_C(u)=1$ and 
$d_C(v)=2$, and a {\em poor couple}  is a $1,2$-couple such that the $2$-vertex is poor, see Figure~\ref{poor-vertices}.  

\iffalse

\vskip .1 true in
\begin{definition}
A {\em $t$-cluster} in an OLD-set $S$ is a connected component in $S$ with $t$ vertices. 
\end{definition}
\vskip .1 true in

We note that $t\ge 2$, as a $1$-cluster is not allowed in an OLD-set. The following definitions assume that the
$t$-cluster $C$ is in an OLD-set $S$.

\vskip .1 true in
\begin{definition}
Let {\em $d_C(v)$} denote the number of neighbors $v$ has in a $t$-cluster $C$.
\end{definition}

\vskip .1 true in
\begin{definition}
A {\em $k$-vertex} is any vertex $v$ in a $t$-cluster $C$ for which $d_C(v)=k$.
\end{definition}

\vskip .1 true in
\begin{definition}
A {\em corner $2$-vertex} in a $t$-cluster $C$ is a $2$-vertex whose two neighbors in $C$ share two neighbors in $G$.
\end{definition}
\vskip .1 true in

\begin{definition}
A {\em $3^+$-vertex} in a $t$-cluster $C$ is a vertex $v$ with $d_C(v)\ge 3$..
\end{definition}
\vskip .1 true in

\begin{definition}
A {\em poor $2$-vertex} is a corner $2$-vertex in a $t$-cluster that has no $3^+$-vertices or corner $2$-vertices as neighbors in $G$.
\end{definition}
\vskip .1 true in

\begin{definition}
A {\em $1,2$-couple} is a pair of adjacent vertices $u,v\in C$ such that $d_C(u)=1$ and $d_C(v)=2$.
\end{definition}
\vskip .1 true in

\begin{definition}
A {\em poor couple} is a $1,2$-couple such that the $2$-vertex is a poor $2$-vertex.
\end{definition}
\vskip .1 true in

\fi

Here are the rules we use to redistribute the weights:

\begin{enumerate}[(R1)]

\item Every vertex gets $\frac{4}{13k}$ charge from each of the $k$ clusters to which it is adjacent.

\item A vertex which gains a charge of $\beta$ from a cluster and is adjacent to $l$ vertices in that 
cluster gains $\frac{\beta}{l}$ from each of those $l$ vertices.

\item In a $t$-cluster $C$ with $t\ge 5$, a $1$-vertex not in a $1,2$-couple receives $7/39$ from its neighbor,  
a poor $2$-vertex not in a poor couple receives $1/39$ from its non-poor neighbors, and a poor couple receives $2/39$ from its neighbor. 

%\item In a $t$-cluster $C$ with $t\ge 5$,  a vertex $u\in C$ with $d_C(u)=2$ gives $2/39$ to its poor $2$-neighbor (corner $2$-vertex with a $1$-neighbor), and $1/39$ to its corner $2$-neighbor. 
\end{enumerate}

\bigskip

Now we show that every vertex will have a weight of at least $\frac{4}{13}$. By doing so, we show that either the vertex has weight at 
least $4/13$, or the vertex belongs to a set of vertices (cluster or couple, for example) that has average weight at least $4/13$. 

% ?? should each vertex not in S have a final weight of AT LEAST 4/13 or exactly 4/13? RIGHT.

Clearly, each vertex not in $S$ has final weight at least $4/13$, by (R1).  For any vertex in $S$, the proof proceeds by considering 
all possible cluster sizes in which the vertex may reside. \\

{\bf Case 1}:  $2$-clusters.  Let the two vertices in this cluster $C$ be $u$ and $v$. Then there are eight vertices not in $S$ and adjacent 
to $u$ or $v$.  Among the eight vertices, two of them are adjacent to both $u$ and $v$, thus at least one of them is adjacent to another cluster, 
then by (R1), these two vertices gain at most $4/13+2/13=6/13$ from $C$.  Each of the three vertices only adjacent to $u$ (and similarly to $v$) 
must be adjacent to another cluster, thus each of them gains at most $2/13$ from $u$.  Therefore, $C$ gives at most $6/13+6\cdot 2/13=18/13$ to 
those vertices.  So the final weights on $u$ and $v$ are $2-18/13=8/13=2\cdot 4/13$.  \\

{\bf Case 2}: $3$-clusters.  The only $3$-cluster is a triangle, for otherwise it is a path in which the two outer vertices would share the same 
neighborhood (the middle vertex) in $S$.

Let the three vertices in the $3$-cluster $C$ be $u, v$ and $w$. Then $C$ has $9$ neighbors not in $S$, three of which are adjacent to two vertices 
in $C$ and six of which are adjacent to exactly one vertex in $C$.  The vertex adjacent to $u,v$ (and similarly to $u,w$ and $w,v$) must be adjacent 
to another cluster, for otherwise it shares the same neighborhood with $w$ in $S$, thus it gains at most $2/13$ from $C$. 

For the two vertices only adjacent to $u$ (and similarly $v$ and $w$), one of them must be adjacent to another cluster, for otherwise they share the 
same neighborhood ($u$) in $S$.  So $C$ gives at most $4/13+2/13=6/13$ to the two vertices and by symmetry, $3\cdot 6/13=18/13$ to the six vertices 
with exactly one neighbor in $C$.  

So the final weights on $u, v, w$ are $3-6/13-18/13=15/13>3\cdot 4/13$. \\

\begin{figure}[h]
\begin{center}
\includegraphics[scale=0.6]{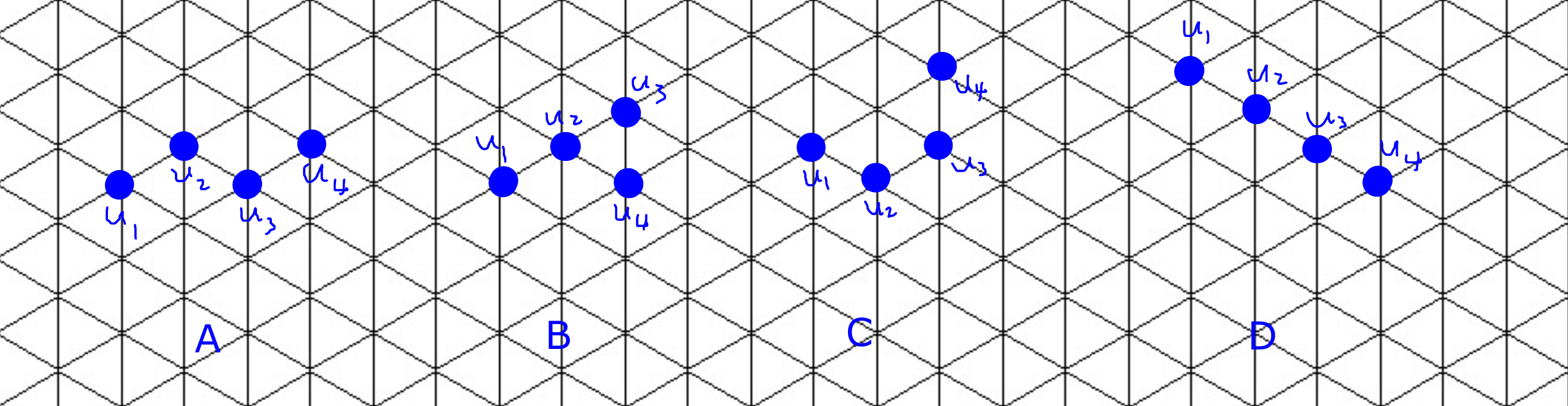}
\caption{All possible $4$-clusters}
\label{4-clusters}
\end{center}
\end{figure}

{\bf Case 3}: $4$-clusters.   The four possible configurations of a $4$-cluster are show in Figure~\ref{4-clusters}. 

%We now consider the four cases individually.  \\

Case 3a:  The $4$-cluster in Figure~\ref{4-clusters}A.  The cluster $C$ has $12$ neighbors outside of $S$, 
and among them, two are adjacent to three vertices in $C$ which gain at most $2\cdot 4/13=8/13$ from $C$ by (R1), two 
are adjacent to two vertices in $C$ which gain at most $2\cdot 4/13=8/13$ from $C$ by (R1), two are adjacent to exactly 
one vertex in $C$ ($u_1$ or $u_3$) which are also adjacent to another cluster thus gaining $2\cdot 2/13=4/13$ from $C$ 
by (R1), and the remaining six vertices are either adjacent to $u_1$ or $u_2$ and four of them are adjacent to 
another cluster, thus gaining $2(4/13+2\cdot 2/13)=16/13$ from $C$ by (R1).  So the final weights on the vertices 
in $C$ are at least $4-8/13-8/13-4/13-16/13=16/13=4\cdot 4/13$. \\ 

Case 3b: The $4$-cluster in Figure~\ref{4-clusters}B.  The vertex adjacent to three vertices in $C$ gains 
at most $4/13$ from $C$,  the vertex adjacent to $u_1, u_2$ and the vertex adjacent to $u_3$ and $u_4$ gain at most 
$2\cdot 4/13=8/13$ from $C$,  the vertex adjacent to $u_2$ and $u_3$ is adjacent to another cluster, thus it gains 
at most $2/13$ from $C$, two of the three vertices adjacent to $u_1$ are adjacent to another cluster, thus the three 
vertices gain at most $4/13+2\cdot 2/13=8/13$ from $C$,  one of the two vertices adjacent to $u_3$ is adjacent to 
another cluster thus the two vertices gain at most $4/13+2/13=6/13$ from $C$, and the vertex adjacent to $u_4$ gains at 
most $4/13$ from $C$. So the final weights on $C$ are at least $4-8/13-2/13-8/13-6/13-4/13=24/13>4\cdot 4/13$. \\

Case 3c: The $4$-cluster in Figure~\ref{4-clusters}C. By (R1), the vertex with four neighbors in $C$ 
gains $4/13$ from $C$, the three vertices adjacent to two vertices in $C$ gain $3\cdot 4/13=12/13$ from $C$,  the 
vertex adjacent to $u_2$ and the vertex adjacent to $u_3$ are adjacent to other clusters, thus they gain at most 
$2\cdot 2/13=4/13$ from $C$, two of the three vertices adjacent to $u_1$ (and similarly to $u_4$) are adjacent 
to other clusters, thus these three vertices gain at most $4/13+2\cdot 2/13=8/13$. So the final weights on $C$ are at 
least $4-4/13-12/13-4/13-2\cdot 8/13=16/13=4\cdot 4/13$. \\

Case 3d: The $4$-cluster in Figure~\ref{4-clusters}D. One of the two vertices adjacent to $u_1$ and 
$u_2$ (and similarly to $u_2,u_3$, or to $u_3,u_4$) must be adjacent to another cluster, thus by (R1) these two vertices 
gain $4/13+2/13=6/13$ from $C$, and two of the three vertices adjacent to $u_1$ (and similarly to $u_4$) are 
adjacent to other neighbors, thus by (R1) these three vertices gain at most $4/13+2\cdot 2/13=8/13$ from $C$. So the 
final weights on $C$ are at least $4-3\cdot 6/13-2\cdot 8/13=18/13>4\cdot 4/13$. \\

{\bf Case 4}: $t$-clusters with $t\ge 5$.  This cluster has initial weight $t$ and, after discharging, should have at least $4t/13$ weight remaining.  
Consequently, the cluster can afford to discharge at most $t-4t/13=9t/13$ and each vertex in the cluster $C$ can afford to discharge $9/13$. 

% ?? Below is a bit unclear here (to me), why can't a vertex w not in S, with neighbors in S, not be in any cluster?
% Do you mean that in previous cases w could be assumed to be adjacent to at least two different clusters?  RIGHT. BECAUSE OF THE SIZE OF CLUSTERS.
% ?? please check that what I changed below is correct.

Note that in the previous cases, when a vertex $w$, outside of $S$, had neighbors in $S$ we could assume that $w$ 
also had neighbors in a different cluster.  But with general clusters of size $t\ge 5$, we cannot assume this.  In fact we must assume 
the worst case, that $C$ must discharge all of its available weight, as the neighbors of $w$ may in fact all lie in the same cluster.
We will show that each poor couple has a final weight of at least $8/13$, and vertices not in poor couples have a final weight of at least $4/13$. 
We consider a few cases, depending on whether a vertex in $C$ has degree one or more. \\

\begin{figure}[h]
\begin{center}
\includegraphics[scale=1.5]{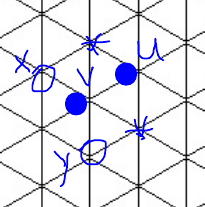}
\caption{Case 4.1}
\label{poor-1-vertex}
\end{center}
\end{figure}

{\bf Case 4.1}: $1$-vertex not in poor couples.  As shown in FIgure~\ref{poor-1-vertex}, let $u$ be a poor $1$-vertex and $v\in C$ be 
its neighbor. By definition, $d_C(v)\ge 3$, so $x$ or $y$ must be in $C$ as well.  Without loss of generality, we may assume that $x\in C$.  

By (R1) and (R2), the vertex not in $S$ and adjacent to $u,v,x$ needs up to $1/3\cdot 4/13$ weight from $u$; the vertex not in $S$ and adjacent 
to $u,v$ needs up to $1/2\cdot 4/13$; $2$ of the $3$ vertices not in $S$ and adjacent to $u$ should have other neighbors in $S$, thus 
the three vertices may need up to $4/13+2\cdot 2/13$ from $u$;  so $u$ may need to give out $4/39+2/13+8/13=9/13+7/39$, but as $u$ 
initially has $9/13$ and gains $7/39$ from $V$ by (R3),  it has enough weight to give out.   \\

% ?? the above is a long sentence that is difficult to parse. Can you break it up into 2 or more sentences?  In particular,
% ?? you say that u has 9/13 intially but gains 7/39 from u...too many u's? $u$ SHOULD BE $v$.

\begin{figure}[h]
\begin{center}
\includegraphics[scale=0.8]{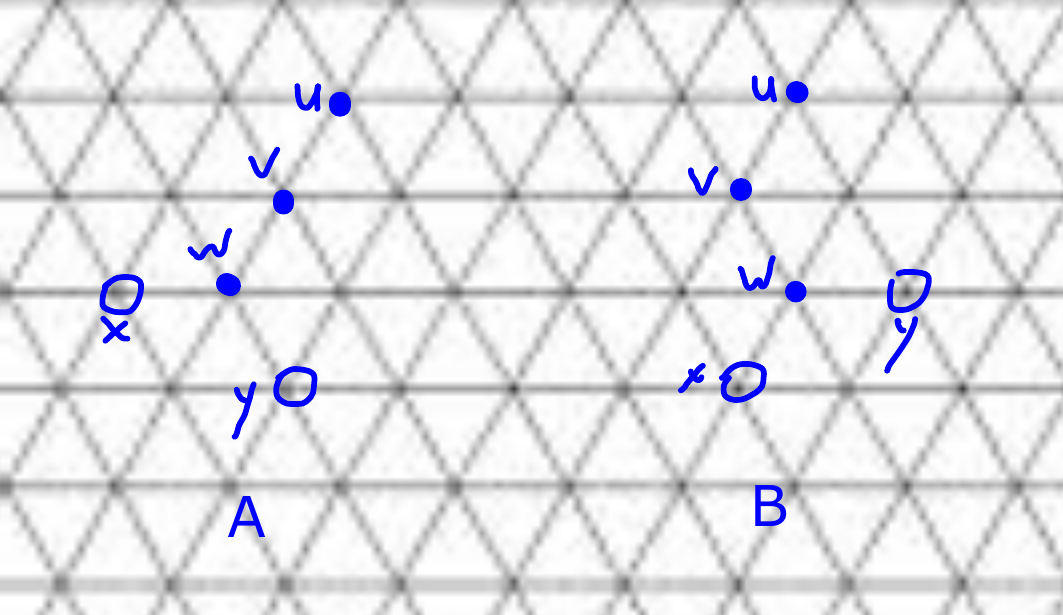}
\caption{Case 4.2}
\label{poor-couple}
\end{center}
\end{figure}

{\bf Case 4.2}: $1,2$-couples.   As shown in the Figure~\ref{poor-couple},  let $u,v$ be a $1,2$-couple in $C$ with $d_C(v)=2, 
d_C(u)=1$, and $w$ is the neighbor of the $1,2$-couple in $C$.    Note that $v$ could be (Figure~\ref{poor-couple}A) or not be 
(Figure~\ref{poor-couple}B)  a corner vertex. 

In Figure~\ref{poor-couple}A, the vertex not in $S$ and adjacent to $v,w,x$ demands at most $1/2\cdot (4/13)$ (if $x\not\in S$)  or $1/3\cdot (4/13)$ (if $x\in S$) from $\{v,u\}$, and similarly for the vertex not in $S$ and adjacent to $v,w,y$, but one of them 
should have a neighbor in $S$ other than $v,w$, we may assume that $v,u$ needs to give $2/13+4/39$ to those two vertices.  One of the 
two vertices not in $S$ and adjacent to $u$ and $v$ should have another neighbor in $S$, so we may assume that $u,v$ gives out up to 
$4/13+2/3\cdot (4/13)$ to those two vertices.  Two of the three vertices not in $S$ and adjacent to $u$ should have other neighbors 
in $S$, thus we may assume that $u,v$ gives out up to $4/13+2/13+2/13$ to those three vertices.  In total,  $u,v$ may give up to 
$(2/13+4/39)+(4/13+8/39)+8/13=18/13$, and they can afford it, as they initially have $2\cdot 9/13$ to give out. 

% ?? hmmm...maybe we should define poor, since you say v is poor below.  Can we just say v is a poor 2-vertex? YES

In Figure~\ref{poor-couple}B,  $v$ is poor $2$-vertex, thus the couple is poor, if neither $x$ nor $y$ is in $C$, and   
by (R1) and (R2), the vertex not in $S$ and adjacent to $v,w$ needs at most $2/13$ from $u,v$,  the vertex not in $S$ and 
adjacent to $u,v,w$ needs at most $2/3\cdot 4/13$ from $u,v$, the vertex not in $S$ and adjacent to $v$ should have another 
neighbor in $S$ and thus demands at most $2/13$ from $v$, the vertex not in $S$ and adjacent to $u,v$ may demand $4/13$ from $u,v$,  
two of the three vertices not in $S$ and adjacent to $u$ should have other neighbors in $S$, thus the three vertices may demand 
up to $4/13+2\cdot 2/13$ from $u,v$.  So, in total, $u,v$ may need to give out $2/13+8/39+2/13+4/13+8/13=56/39=18/13+2/39$.  As 
$u,v$ initially have $18/13$ and gain $2/39$ more from $w$ by (R3), they have enough weight to give out. 

If $x$ is in $C$ in Figure~\ref{poor-couple}B, then the vertex not in $C$ but adjacent to $v,w,x$ (if $x\in C$) gets 
$1/3\cdot 4/13=4/39$ from $v$, instead of $1/2\cdot 4/13$ as in the poor case. So, the $1,2$-couple may need to give out 
$2/13+8/39+4/39+4/13+8/13=18/13$, which they can afford.   Similarly, if $y$ is in $C$, then the vertex adjacent to 
$u,v,w,y$ may need $1/2\cdot 4/13$ from $u$ and $v$, instead of $2/3\cdot 4/13$ as in the poor case, so the $1,2$-couple 
may need to give out $2/13+2/13+2/13+4/13+8/13=18/13$, which they can afford.  \\

%{\bf Case 4.3}: $1,2$-couples that are not poor.  

{\bf Case 4.3}: $2$-vertex not in $1,2$-couples.   Let $u$ be such a vertex and $v,w$ be the neighbors of $u$. By definition, 
$u$ has no neighbor of degree $1$ in $C$, so $d_C(v), d_C(w)\ge 2$.  

\begin{figure}[h]
\begin{center}
\includegraphics[scale=0.7]{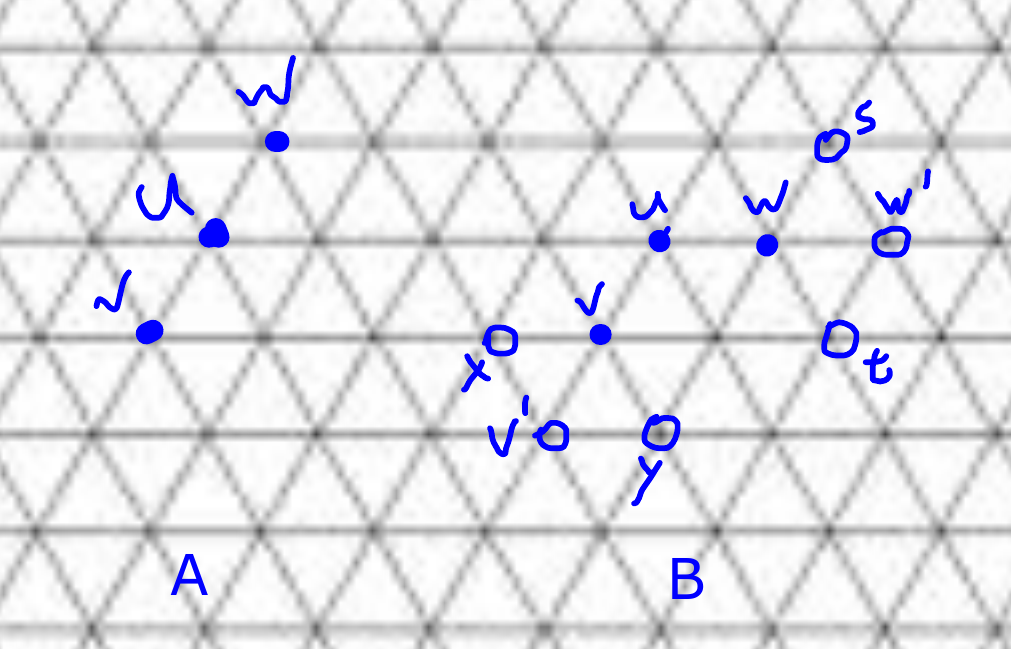}
\caption{Case 4.3}
\label{2-vertex}
\end{center}
\end{figure}

If $v,w$ have no common neighbor (see Figure~\ref{2-vertex} A), then $u$ is not a $2$-corner vertex, and in this case, one of the two 
neighbors $u,v$ share has a third neighbor in $S$, and one of the two neighbors $u,w$ share also has a third neighbor in $S$.  Thus, 
$u$ needs to give out at most $1/2\cdot 4/13\cdot 2+1/3\cdot 4/13\cdot 2=5/3\cdot 4/13$ to the four neighbors outside of $S$, and 
may need to give $2\cdot 2/39$ to $v$ and $w$ (if they are in poor couples).  In total $u$ may need to give out up to $8/13$,  which 
it can afford as it has $9/13$ to start with.   

Next, we assume that $u$ is a corner $2$-vertex, that is, $v$ and $w$ share a common neighbor, see Figure~\ref{2-vertex} B. 

First, let $u$ be {\em not poor}, that is, $v$ or $w$ is a $3^+$-neighbor or a corner $2$-vertex.   
We may assume that $v$ is a $3^+$-vertex or a corner $2$-vertex, thus $x$ or $y$ is in $C$.  
Note that $v$ and $w$ cannot be poor, as they have $u$ as a neighbor and $u$ is a corner $2$-vertex. 

If $x$ is in $C$, by (R1) and (R2), $u$ needs to give out at most $4/13+2/13+1/3\cdot 4/13\cdot 2=2/3<9/13$ to the neighbors 
not in $S$ and by (R3) gives nothing to $v$.  Similarly if $y$ is in $C$, then $u$ gives out at most $4/13+2/13+2/13+1/4\cdot 4/13=9/13$. 
Thus, $u$ can afford both cases since it has $9/13$ initially.  
 
Now let $u$ be {\em poor}.  Then both $v$ and $w$ are $2$-vertices and not corner $2$-vertices.  
If $C$ is not a $5$-cluster, then $v$ or $w$ is not in a
$1,2$-couple, thus not poor, and so by (R3) $u$ gains $1/39$ from it.  Therefore, $u$ need to give out at most 
$4/13+2\cdot 2/13+1/3\cdot 4/13=9/13+1/39$, but can afford to do so, as it gains $1/39$ and has $9/13$ initially.  

If $C$ is a $5$-cluster, then let the other neighbors of $v,w$ be $v',w'$ respectively.  Then by (R1) and (R2), 
$C$ needs to give at most $4/13$ to the vertex not in $S$ but adjacent to only $u$, or only $u,v,w$, or only 
$u,v$, or only $u,w$.  One of the two common neighbors of $w$ and $w'$ should have a neighbor in $S-C$, thus 
$C$ gives out at most $1/2\cdot 4/13+4/13$ to those two neighbors, and similarly gives at most $1/2\cdot 4/13+4/13$ to 
the two common neighbors of $v$ and $v'$.  Two of the $3$ neighbors of $w'$ (and similarly of $v'$)
should have another neighbor in $S$, thus $C$ gives out at most $(4/13+2\cdot 1/2\cdot 4/13)\cdot 2=16/13$ to 
the $6$ neighbors of $w'$ and $v'$.  Therefore, $C$ gives out at most $4/13\cdot 4+2\cdot (1/2\cdot 4/13+4/13)+16/13=44/13$, 
but $5-5\cdot 4/13=45/13>44/13$, so $C$ can afford to do so.  

% ?? in the above, I don't follow what you mean when you write "Two of the three neighbors only to w'"
% ?? do you mean "2 of the 3 neighbors of w'"  ? YES === changed "only to" to "of" rkk

%As shown in Figure~\ref{2-vertex},  a $2$-vertex in $S$ but not in poor couples has four neighbors not in $S$ and may be 
%adjacent to two poor couples (shown in the figure).  By (R1) and (R2) and (R3), one can see that in Figure~\ref{2-vertex} D,  
%$u$ needs to give out most weights, which are $4/13+2/13+2/13+4/39=9/13+1/39$ to vertices not in $S$ and $2\cdot 2/39$ to poor 
%couples,  \red{the poor couples in D does not need charges. the 5-cluster is in trouble!!}

\begin{figure}[h]
\begin{center}
\includegraphics[scale=1.2]{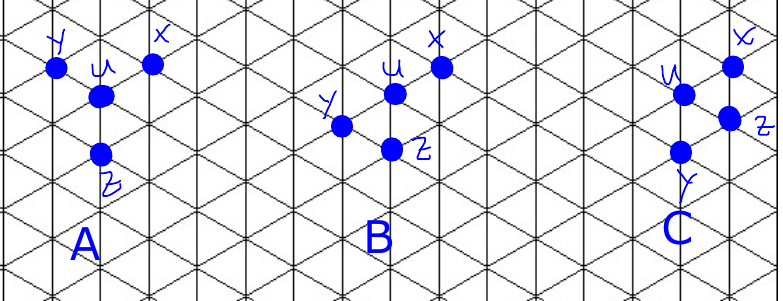}
\caption{Case 4.4}
\label{3-vertex}
\end{center}
\end{figure}

{\bf Case 4.4}: $3^+$-vertices.  As shown in Figure~\ref{3-vertex}, we assume that $u$ with $d_C(u)=3$ 
has three neighbors $x,y,z$ in $C$.   Note that at most one of $x,y,z$ has degree $1$ in $C$, but all 
of them could be in {\em poor couples}.  Thus by (R1), (R2) and (R3), we see that in Figure~\ref{3-vertex}A, $u$ gives 
out at most $4/13\cdot 1/3\cdot 3+7/39+2\cdot 2/39=23/39<9/13$; in Figure~\ref{3-vertex}B,  $u$ gives 
out at most $2/13\cdot 2+4/13\cdot 1/3+7/39=23/39<9/13$;  and in Figure~\ref{3-vertex}C,  $u$ gives out 
at most $4/13+2/13\cdot 2=8/13<9/13$.  

Thus every vertex in $S$ ends up with a weight at least $4/13$ as well.  

\section{Future Research}

In this paper, we give the exact OLD-density of the infinite triangular grid.   
The OLD-density problem of other infinite grids, especially non-regular ones, 
could also be interesting.  For example, what is the OLD-density of infinite 
triangular-hexagonal grids?  In the study of identifying codes, $r$-identifying codes, 
meaning a code-vertex can cover vertices within distance $r$, have been considered, see ~\cite{MS10}.  
In the context of OLD-sets, one may also analogously study an $r$-OLD-set.  

\section*{Acknowledgement}
We would like to thank Pete Slater for his valuable comments and suggested references.

\end{document}